\documentclass[12pt,a4paper]{article}

\usepackage{amsmath}
\usepackage{amssymb}
\usepackage{amsthm}
\usepackage{amsfonts}
\usepackage{graphicx}
\usepackage[all]{xy}
\usepackage{verbatim}
\usepackage{parskip}
\usepackage{csquotes}
\usepackage{tkz-graph}
\usepackage[left=2cm,right=2cm,top=3cm,bottom=2.5cm]{geometry}
\usepackage{soul}

\usepackage{hyperref}
\usepackage{ulem}

\newcommand{\norm}[2]{\left\| #1 \right\|_{#2}}
\newcommand{\M}{\mathcal{M}}
\newcommand{\Om}{\Omega}
\newcommand{\R}{\mathbb{R}}
\newcommand{\Z}{\mathbb{Z}}
\newcommand{\C}{\mathbb{C}}
\newcommand{\N}{\mathbb{N}}

\theoremstyle{definition}
\newtheorem{theorem}{Theorem}[section]
\newtheorem{definition}[theorem]{Definition}

\newtheorem{remark}[theorem]{Remark}
\newtheorem{proposition}[theorem]{Proposition}

\newtheorem{corollary}{Corollary}[section]
\newtheorem{lemma}[theorem]{Lemma}

\begin{document}
%%%%%%%% Title Block %%%%%%%%%%%%%%
\thispagestyle{empty}
\begin{center}
    {\LARGE On Bridging Analyticity and Sparseness in Hyperdissipative Navier--Stokes Systems \par}

    \vspace{0.6cm}

    {\large
    Moses Patson Phiri \par}
    \vspace{0.2cm}

    {\normalsize
    Department of Mathematics,
    University of Alabama at Birmingham
    \par}

    \vspace{0.4cm}

    {\normalsize \today\par}
\end{center}

\vspace{0.4cm}
%\hrule
\vspace{0.8cm}

\begin{abstract}
We study the three-dimensional hyperdissipative Navier–Stokes system in the super-critical regime (below the Lions threshold). Leveraging the analyticity–sparseness gap in the case where the ratios of the lower-order derivatives with respect to the higher-order derivatives are bounded by constants, we introduce time-weighted bridge inequalities quantifying these ratios across all derivative levels. On one hand, the time weights are chosen in a way that the resulting lower bound on the radius of spatial analyticity still dominates the scale of sparseness of the super-level sets indicating that the flow entered the dissipation regime. On the other hand, they provide a less restrictive environment for the ratios giving the construction more flexibility. As an illustration of its utility, it is shown that -- together with a harmonic-measure contraction on one-dimensional sparse sets -- this mechanism rules out a more general class of the blow-up scenarios compared to the previous work by Farhat and Gruji\'c (essentially a perturbation of the multidimensional geometric series).
\end{abstract}


 \newpage   
\section{Introduction}  \label{section1}
The three-dimensional Hyperdissipative  Navier-Stokes  system in \(\R^3 \times (0, T)\) is given by
\begin{align}
\label{hdnse}
\begin{cases}
\partial_t u + (-\Delta)^\beta u + (u\cdot \nabla) u + \nabla p = 0,\\
\nabla \cdot u =0,\\
u(\cdot , 0) = u_{0}(\cdot).
\end{cases}
\end{align}
The nonlinearity can be written in the equivalent divergence form 
 \[(u\cdot \nabla) u := \sum_{j = 1}^{3} u_j \partial_j u_i =  \sum_{j = 1}^{3}\partial_{j}(u_{i}u_{j})  = \nabla \cdot \left(u \otimes u \right).\]
The exponent \(\beta > 1\) controls the order of dissipation, the operator \((-\Delta)^\beta\) acts in Fourier space as multiplication by \(|\xi|^{2\beta}\), the vector field \(u: \R^3 \times (0,T) \to \R^3\) is the velocity of the fluid, and the scalar field \(p: \R^3 \times (0,T) \to \R^3\) is the pressure. For simplicity, the kinematic viscosity \(\nu\) is normalized to \(1\), and the external force to zero.

Under the rescaling
\[u_\lambda(x,t)=\lambda^{2\beta-1}u(\lambda x,\lambda^{2\beta}t),  \qquad
p_\lambda(x,t)=\lambda^{4\beta-2}p(\lambda x,\lambda^{2\beta}t),\]
the system (\ref{hdnse}) is invariant. It is critical in three dimensions when \(\beta=\tfrac{5}{4}\) (i.e., the energy norm is invariant), subcritical for \(\beta > \tfrac{5}{4}\), and supercritical for \(1\leq \beta<\tfrac{5}{4}\).

Global well-posedness for the hyperdissipative Navier-Stokes equations is classical 
for exponents $\beta \geq \frac{d+2}{4}$, a result originating with 
J.L. Lions \cite{Lions1969}. In three dimensions, this corresponds to $\beta \geq 5/4$. Building on this boundary, Tao \cite{Tao2009} proved that global regularity persists even under logarithmically supercritical weakening of the critical dissipation, replacing
\[(-\Delta)^{\tfrac{5}{4}} \qquad \text{by} \qquad \tfrac{(-\Delta)^{\tfrac{5}{4}}}{(\log(2-\Delta))^{\tfrac{1}{2}}}. \]
Thus, while the Lions exponent \(\tfrac{5}{4}\) marks the energy–critical threshold in 3D, even a half-power logarithmic deficit remains controllable, as rigorously confirmed by the logarithmically supercritical global-regularity result of Barbato–Morandin–Romito \cite{BarbatoMorandinRomito}. Furthermore, Colombo--Haffter proved that, for any fixed divergence-free datum $u_0 \in H^\delta$, there exists an \(\varepsilon > 0\) which depends on \(\delta\) and the initial datum, such that one retains global smooth solutions for all fractional exponents $\beta \in (5/4 - \varepsilon, 5/4]$ \cite{Colombo-haffter}. In the same hyperdissipative regime \(1< \beta <5/4\), Katz–Pavlović \cite{Katz-Pavlovic} obtained a Caffarelli–Kohn–Nirenberg-type partial regularity result, showing that at the first blow-up time the spatial singular set has Hausdorff dimension at most \(5-4\beta\).

In the supercritical regime \(1 \leq \beta<\tfrac{5}{4}\), unconditional global regularity remains open. Recent progress leverages structure suggested by spatial intermittency of turbulent flows, geometric sparseness of active sets, to explore conditions under which the scale of sparseness falls below the scale of the radius of spatial analyticity (which is a natural dissipation scale) depleting the nonlinearity at dynamically relevant scales. In recent work, Gruji\'c and Xu utilized their mathematical framework based on the scale of sparseness of the super-level sets of higher-order derivatives to show that as soon as $\beta>1$ no blow-up can occur in `turbulent scenario' \cite{grujic}. The mathematical condition is that the higher-order derivatives dominate lower-order derivatives in a suitable dynamic sense and across a finite range of indices (the closer $\beta$ to 1 the wider the range); these derivative ratios can be viewed as multi-level, local versions of the classical Taylor or Kraichnan scales in turbulence phenomenology (gradients dominating the fields), hence `turbulent'. In particular, this work ruled out asymptotically self-similar blow-up, a prime candidate for singularity formation, as soon as $\beta>1$. The proof is based on carefully tracking ordered dynamics of the portions of the `chain of derivatives', ascending (higher-order derivatives dominate, consistent with the turbulent buildup) vs. descending (lower-order derivatives dominate, consistent with the turbulent decay). The utility of the ascending scenario is that it replaces the use of classical Gagliardo-Nirenberg interpolation at higher differential level which is precisely what yields a scaling gain compared to the classical methods; this was termed `dynamic interpolation'.

Since the arguments in the general case presented in \cite{grujic} are quite technical and since the ascending scenario is the one more relevant in the context of singularity formation, Farhat and Gruji\'c in \cite{farhat} illustrated the mechanism behind the scaling gain in the ascending scenario in the case where the singularity buildup can be described as a perturbation of a multidimensional Taylor series for the velocity, the radius of analyticity shrinking to 0 as the flow approaches the singular time (there is no assumption on the rate at which the analyticity radius shrinks to zero, i.e., on the strength of the singularity). Ruling out this type of the blow-up does follow from the general theory presented in \cite{grujic}; the main point in \cite{farhat} was to have a self-contained and transparent presentation in a simple setting.

The present paper continues this line of work in the supercritical regime. We show that the ascending building-block mechanism still rules out blow-up even when the relevant derivative ratios are allowed to grow with respect to time as the possible singular time is approached. Thus, the constant bounds used in the predecessor construction are replaced by time weighted bounds (\ref{time-weight1}), while the key analyticity--sparseness comparison is preserved. We now state the main result. 

\begin{theorem} \label{main-theorem}
Let \(\beta > 1\), \(u_0 \in L^{\infty} \cap L^2\), denote by \(T^*\) the first possible singular time, and by \(x^*\) an isolated spatial singularity at \(T^*\). Then, for any \(t \in (0,T^*)\), the solution \(u(t)\) is spatially analytic \cite{Guberovic2010}, and each velocity component \(u^i\) admits the expansion
\[
u^i(x,t)=\sum_{k=0}^{\infty}\sum_{|\alpha|=k} c^i_{\alpha,k}(t)(x-x^*)^\alpha
\]
in a neighborhood of \(x^*\), where \(\alpha\) is the multi-index \(\alpha=(\alpha_1,\alpha_2,\alpha_3)\). %Similarly to \cite{farhat}, 
We make the following assumptions on the form of a singularity buildup.

       \begin{enumerate}
        \item[(i)] suppose that the Taylor coefficients have the following form: there exists constants \(\varepsilon > 0,\) \(C > 1\), and \(\gamma >0\) to be determined in the proof such that for all \(t \in (T^* - \varepsilon, T^*)\), \(c_{\alpha, k}^{i}(t) := \delta_{\alpha, k}^{i}(t) \dfrac{1}{\rho(t)^k}\), 
        where \(\rho > 0, \rho \to 0\) as \(t \to T^{*}\), and for \(k \not=0\), we have 
        \begin{align}\label{modulation}
            \dfrac{1}{C^k (T^{*} - t)^{-\frac{\gamma k }{\log k} } } \le \delta_{\alpha, k}^{i} \le C^k (T^{*} - t)^{-\frac{\gamma k}{\log k} }
        \end{align}

       with the additional assumption that,  for every relevant \(i, \alpha,\) and \(k,\) the coefficient function
       \(t \mapsto c_{\alpha, k}^{i}(t)\) is increasing on \((T^{*}-\varepsilon, T^{*}).\)
        \item[(ii)] \label{cod2} suppose that the blow-up is focused at \(x^{*}\):  for all \(t \in (T^*-\varepsilon, T^*)\) we have
        \[\|D^{(k)}u(t) \|_{\infty} = |D^{(k)}u(x^{*}, t)|, \qquad \text{for all } k \in \N.\]
    \end{enumerate}
  Then \(T^*\) is not a singular time, and the solution \(u(x,t)\) extends analytically beyond \(T^{*}.\)
\end{theorem}

\begin{remark}
    The assumption that the coefficient functions \(c_{\alpha, k}^{i}\) are increasing in time is an additional structural hypothesis on the singularity buildup. It is not a consequence of the representation 
    \(c_{\alpha, k}^{i}(t) = \delta_{\alpha, k}^{i}(t) {\rho(t)^{-k}} \), or of the condition that \(\rho(t) \to 0\) as \(t \to T^{*}\), since no rate of decay is imposed on \(\rho.\)

    This hypothesis enters only at the final stage of the proof. Together with the focused-buildup condition, it provides the one-sided temporal growth of the relevant level-\(k\) derivative norm and hence the escape-time property needed to contradict the strict contraction obtained from the harmonic-measure maximum principle. We emphasize that this assumption is restrictive and may exclude more intermittent singularity scenarios in which different frequencies or derivative levels become dominant at different times. Weakening this assumption while preserving the escape-time mechanism is an interesting direction for future work.
\end{remark}

\begin{remark}
The focused blow-up hypothesis (ii) (as in \cite{farhat}) should be viewed as a strong localization assumption, presented here as a model scenario for highly concentrated singular behavior. Its role is to test whether a blow-up mechanism in which the extremizing derivatives remain centered at a single spatial point can persist up to \(T^*\). The key difference between our hypothesis (i) and the one in \cite{farhat} is that the modulation inequality \eqref{modulation} is now time-dependent.
\end{remark}

We now explain the meaning of the assumptions in the theorem and how the result fits into the preceding analyticity--sparseness framework. The main building blocks in the theory presented in \cite{grujic} are the lower-to-higher ascending portions and the higher-to-lower descending portions of the chain of derivative ratios. In the original construction, these ratios are bounded by constants (see (\ref{bddconstants})). The ascending case is the one relevant to singularity formation: higher-order derivatives dominate lower-order derivatives, reflecting the buildup of small scales.

In the ascending case, this structure replaces the use of the classical Gagliardo--Nirenberg interpolation when estimating the lower-order terms arising from the Leibniz expansion of the Navier--Stokes nonlinearity in the complexified Duhamel formula for higher-order derivatives. This replacement is essential: it produces the scaling gain that improves the lower bound on the radius of spatial analyticity as soon as \(\beta>1\). Consequently, at sufficiently high differential levels \(k\), the resulting
analyticity-radius scale dominates the corresponding level \(k\) sparseness scale, thereby pushing the flow back into the dissipation-dominated regime.

The main goal of the present paper is to generalize the ascending building-block design introduced by Grujić. Instead of requiring the relevant derivative ratios to be bounded by constants, that is
\begin{align}\label{bddconstants}
    \frac{\norm{D^ju(t)}{\infty}^{\frac{1}{j+1}}}{\norm{D^ku(t)}{\infty}^{\tfrac{1}{k+1}}}
     \le  \dfrac{\left(j!\right)^{\frac{1}{j+1}}}{\left(k!\right)^{\tfrac{1}{k+1}}} , \qquad \ell \le j \le k, 
\end{align}
we allow them to be controlled by algebraic time-weights that may grow without bound as the possible singular time is approached. The purpose of the time-weighted bridge inequalities is to preserve the scaling gain described above under this more flexible hypothesis. More precisely, the time-weighted bridge inequalities have the form

\begin{align}\label{time-weight1}
    \frac{\norm{D^ju(t)}{\infty}^{\frac{1}{j+1}}}{\norm{D^ku(t)}{\infty}^{\tfrac{1}{k+1}}}
     \le (T-t)^{-\tfrac{1}{2\beta} \left(\tfrac{\mu + j}{j+ 1} - \tfrac{\mu + k}{k + 1}\right)} \dfrac{\left(j!\right)^{\frac{1}{j+1}}}{\left(k!\right)^{\tfrac{1}{k+1}}}, 
\end{align}
for all \(  \ell \le j \le k \) and \(t\in(0,T)\), where \(\mu>1\) is a suitable parameter. Since 
\[ \frac{\mu+j}{j+1}-\frac{\mu+k}{k+1}
=(\mu-1)\!\left(\frac{1}{j+1}-\frac{1}{k+1}\right) \geq 0,
 \]
for any  \(\mu>1\) and \(j \leq k,\) the time factor in (\ref{time-weight1}) allows for much stronger fluctuations in the derivative ratios as \(t\uparrow T\). This time-dependent flexibility is a key distinction between the construction in Farhat-Grujic \cite{farhat} and the one developed in this paper.

In particular, the construction in Farhat-Grujic yields the analyticity-radius estimate 
 \[
\rho_k \;\gtrsim\; \|D^k u\|_{\infty }^{-\,\frac{1}{2\beta-1}\frac{1}{k+1}},
\] 
whereas the present time-weighted construction yields
\[
\rho_k \;\gtrsim\; \|D^k u\|_{\infty }^{-a_k},
\qquad
a_k=
\frac{1}{2\beta-1-(\mu-1)\frac{2k}{k+1}}\frac{1}{k+1}.
\]  

Although the time-dependent control produces a weaker lower bound on the radius of spatial analyticity than in the constant-coefficient construction, the bound remains strong enough, for sufficiently large differential level (k), to dominate the corresponding level (k) sparseness scale. This analyticity--sparseness domination is the key mechanism that prevents the formation of a singularity.

As an application, we show that these more flexible building blocks (\ref{time-weight1}) allow the potential singularity buildup described in \cite{farhat} to be made more general. More precisely, the buildup considered in \cite{farhat} was modeled as a perturbation of a multidimensional geometric series in which the perturbation was relatively rigid: the modulation of the geometric-series ratio was time-independent. The time-weighted bridge inequalities  introduced here allow this modulation to depend on time.

\begin{remark}
One should note that in \cite{AlbrittonBradshaw},  Albritton and Bradshaw introduced an alternative way to study sparseness both in the physical and frequency scales. This was presented in the case of the Navier-Stokes system and on the velocity level; it would be interesting to see if their setting could be used to obtain regularity results analogous to the ones presented in \cite{grujic} and -- in particular -- to show that a suitably formulated analogue of the ascending scenario would lead to regularity for the hyperdissipative Navier-Stokes system in the super-critical regime.
\end{remark}

This paper is organized as follows: Section \ref{prelim} collects preliminaries on analytic classes, the 1D \(\delta-\)sparseness framework, and harmonic-measure tools. In section \ref{Mainsection}, we reintroduce the above time-weighted bridge inequalities and present the two main results: the improved (compared to the classical theory) lower bound on the radius of spatial analyticity in the ascending scenario delineated by the time-weighted bridge inequalities and the application to Farhat-Gruji\'c blow-up scenario.

\section{Preliminaries}\label{prelim}

In this section, we collect the notions and baseline facts on \textit{sparseness of regions of intense fluid activity} that we use throughout the paper. The quantitative framework for sparse superlevel sets in 3D was initiated by Grujić \cite{grujicz} and further refined in subsequent works, including the reformulation in \cite{grujic2}, as well as developments by Farhat et al. \cite{farhat2}; see also related refinements by Bradshaw and coauthors \cite{Bradshaw}. Our presentation follows these references while tailoring the statements to the hyperdissipative setting.

\paragraph{Notation and conventions.}
We write \(D^k:=\nabla^k\) for the \(k\)th spatial derivatives, and \(\|\,\cdot\,\|_p\) for \(L^p(\R^3)\) norms; in particular, \(\|\,\cdot\,\|:=\|\,\cdot\,\|_\infty\). We use \(A\lesssim B\) to mean \(A\le CB\) for a constant \(C>0\) independent of \(k\) and \(t\). 
Given a measurable set \(E\subset\R^3\), \(|E|\) denotes its Lebesgue measure.

Firstly, we recall what is meant by local sparseness on a given scale (Gruji\'c, \cite{grujicz}).
\begin{definition}[1D $\delta$-sparse]
Let $S$ be an open set in $\mathbb{R}^3$, $\delta \in (0,1)$, $x_0 \in \mathbb{R}^3$, and $r \in (0,\infty)$. 
We say that $S$ is \textit{$1$-dimensional $\delta$-sparse around $x_0$ at scale $r$} if there exists a unit vector $\nu \in S^2$ such that
\[
\frac{\lvert S \cap (x_0 - r\nu,\, x_0 + r\nu)\rvert}{2r} \le \delta,
\]
where $\lvert\cdot\rvert$ denotes one-dimensional Lebesgue measure (length) along the segment.
\end{definition}

\begin{remark}
In applications, \(S\) will be a super-level set of a component of a derivative magnitude (e.g., \(S=\{x:\ |D^k u(x,t)|>\lambda\}\)). One-dimensional sparseness captures quantitative smallness along at least one direction through a given point and underlies harmonic-measure arguments used to propagate uniform bounds; see \cite{grujic, grujic2, farhat2}.
\end{remark}

\begin{definition}[Escape Time]
    Let \(u \in C([0, T^*), L^\infty)\), where \(T^*\) is the first possible blow-up time.  A time \(t \in (0, T^*)\) is called an scape time if 
    \[\norm{u(\tau)}{} > \norm{u(t)}{} \quad \text{ for every } \quad \tau \in (t, T^*).\]
    The local-in-time well-posedness of the \(L^{\infty}\)-norm implies that for any level there exists a unique escape time.
\end{definition}

Let \(\Omega\subset\mathbb{C}\) be open and \(K\subset\partial\Omega\) closed, and suppose
\(\Omega\) is sufficiently regular. Then there exists a unique
function \(\omega:= \omega(\,\cdot\,, \Omega, K)\), harmonic in \(\Omega\), with
\(\omega=1\) on \(K\) and \(\omega=0\) on \(\partial\Omega\setminus K\).
For \(z\in\Omega\), the value \(\omega(z, \Omega, K)\) is called the harmonic measure of \(K\)
with respect to \(\Omega\) at \(z\). Below we record several basic properties of harmonic measure.

\begin{proposition}[Ransford \cite{Ransford}]
Let \(\Omega\) be an open, connected set in \(\mathbb{C}\) such that its boundary has nonzero Hausdorff dimension, and let \(K\) be a Borel subset of the boundary. Suppose that \(u\) is a subharmonic function on \(\Omega\) satisfying
\[
u(z) \leq M, \quad \text{for } z \in \Omega
\]
and
\[
\limsup_{z \to \xi} u(z) \leq m, \quad \text{for } \xi \in K.
\]
Then
\[
u(z) \leq m h(z,\Omega,K) + M(1 - h(z,\Omega,K)), \quad \text{for } z \in \Omega.
\]
\end{proposition}

\begin{proposition}\label{conformal}
    The harmonic measure is invariant with respect to conformal mappings.
\end{proposition}

\begin{proposition}
    Let \(f\) be a holomorphic function on an open set \(S\) in \(\C,\) then \(\log|f|\) is subharmonic on \(S\).
\end{proposition}

We recall the following theorem of Gruji\'c and Xu \cite{grujic}, which identifies the a priori sparseness scale for the super-level sets of the \(k\)-th order derivatives of the velocity field.

\begin{theorem}[Gruji\'c--Xu \cite{grujic}, sparseness-scale estimate] \label{scalebound}
	Let \(u\) be a Leray solution on \((0,T^*)\) with
	\[
	u\in C((0,T^*);L^\infty(\mathbb{R}^3)).
	\]
	Let \(\zeta\in\mathbb{N}^3\) satisfy \(|\zeta|=k\). For each \(t\in(0,T^*)\), then the super-level sets
	\[
	S_{\zeta,\lambda}^{i,\pm}(t)
	:=
	\left\{
	x\in\mathbb{R}^3 : (D^\zeta u)_i^\pm(x,t)>\lambda \|D^\zeta u(t)\|_{L^\infty}
	\right\},
	\qquad 1\le i\le 3
	\]
	 are \(3D\) \(\delta\)-sparse around any spatial point \(x_0\) at scale
	\[
	r^*(t)=c(\|u_0\|_{L^2})\,\|D^\zeta u(t)\|_{L^\infty}^{-\frac{2}{2k+3}},
	\]
	whenever \(r^*(t)\in(0,1]\).
\end{theorem}

\begin{lemma}[Gagliardo-Nirenberg, \cite{Gagliardo}\cite{Nirenberg}] Suppose \(p, q, r \in [1, \infty], \theta \in \mathbb{R}\), and \(m, j, d \in \mathbb{N}\) satisfy 
\[ \dfrac{1}{p}  = \dfrac{j}{d} + \theta \left( \dfrac{1}{r} - \dfrac{m}{d}\right) + \dfrac{1-\theta}{q}, \quad \dfrac{j}{m} \leq \theta \leq 1.\]
Then, there exists a constant \(C = C(m, d, j, q, r, \theta)\) such that
\begin{align}
    ||D^jf||_{L^p} \leq C ||D^mf||_{L^r}^{\theta} ||f||_{L^q}^{1-\theta}
\end{align}
\end{lemma}

The following theorem of Grujić and Xu \cite{grujic} provides a general mechanism: when a sufficiently large portion of the derivative chain \( \{D^{(j)}u\}_{j\leq k} \) exhibits ascending behavior--meaning that higher-order derivatives dominate lower--order ones--the local existence time is extended, yielding improved estimates for the analyticity radius.

\begin{theorem}[Grujic and Xu \cite{grujic}, Fahart and Grujic \cite{farhat}] \label{theorem1}
Let $\beta>1$, $u_0\in L^2$, and $D^{i}u_0\in L^\infty$ for $0\le i\le k$, and suppose that
\begin{equation}\label{notime}
  \norm{D^{m}u_0}{}^{\tfrac{1}{m+1}}
  \;\leq\;
  \M_{m,k}\,\norm{D^{k}u_0}{}^{\tfrac{1}{k+1}}
  \qquad \forall\, \ell \le m \le k .
\end{equation}
Assume further that the constants $\{\M_{j,k}\}$ and the indices $\ell$ and $k$ satisfy
\begin{equation}\label{cond1}
  \sum_{\ell \le i \le j-\ell}
  \binom{j}{i}\,
  \M_{i,k}^{\,i+1}\,
  \M_{j-i,k}^{\,j-i+1}
  \;\lesssim\; \phi(j,k),
  \qquad \forall\,\, 2\ell \le j \le k ,
\end{equation}
and
\begin{equation}\label{cond2}
  \norm{u_0}{2} \!
  \sum_{0 \le i \le \ell}
  \binom{j}{i}
  \M_{\ell,k}^{\,\frac{(\ell+1)(i+\frac{3}{2})}{\ell+\frac{3}{2}}}
  \M_{j-i,k}^{\,j-i+1}
  \Big((k!)^{\frac{1}{k+1}}\norm{u_0}{}\Big)^{\frac{\left(\frac{3}{2}-1\right)(\ell-i)}{\ell+\frac{3}{2}}}
  \lesssim \psi(j,k),
\end{equation}
 for all  \( 2\ell \le j \le k, \) for some functions \(\phi\) and \(\psi\). If
\begin{equation}\label{eq:3.4}
  T_j \;\lesssim\;
  \big(\phi(j,k)+\psi(j,k)\big)^{-\tfrac{2\beta}{2\beta-1}}\,
  \norm{D^{k}u_0}{}^{-\tfrac{2\beta}{(2\beta-1)(k+1)}},
\end{equation}
then for any $\ell \le j \le k$ the complexified solution satisfies the upper bound
\begin{equation}\label{eq:3.5}
  \sup_{t\in(0,T_j)} \;\sup_{y\in \mathcal{D}_t} \norm{D^{j}u(\cdot,y,t)}
  \;+\;
  \sup_{t\in(0,T_j)} \;\sup_{y\in \mathcal{D}_t} \norm{D^{j}v(\cdot,y,t)}
  \;\lesssim\;
  \norm{D^{j}u_0}{} + \norm{D^{k}u_0}{}^{\frac{j+1}{k+1}},
\end{equation}
where \[\mathcal{D}_{t} = : \{x + iy \in \mathbb{C}^3 : |y| \le ct^{\frac{1}{2\beta}}\}.\]
\end{theorem}

The coefficient \(\mathcal{M}_{m,k}\) appearing in Theorem \ref{theorem1} is the one defined in \eqref{Mjk-def}.

\section{Main Estimates and Proof of the Main Theorem}
\label{Mainsection}
We now develop the estimates required to prove Theorem \ref{main-theorem}, stated in the Introduction. We start by converting the time-weighted derivative-ratio bounds into a quantitative lower bound on the radius of spatial analyticity, then compare this analyticity scale with the corresponding sparseness scale and complete the proof by applying the harmonic-measure maximum principle.

As outlined in the Introduction, we modify the ascending-chain construction by allowing the ratios between derivative levels to be controlled by algebraic time weights. Accordingly, for \(0\leq j\leq k\), we define the time-weighted coefficient
\begin{equation}\label{Mjk-def}
\mathcal{M}_{j,k}
:=
\bigl(\tau_k\bigr)^{-\frac{1}{2\beta}
\left(\frac{\mu+j}{j+1}-\frac{\mu+k}{k+1}\right)}
\frac{(j!)^{\frac{1}{j+1}}}{(k!)^{\frac{1}{k+1}}}.
\end{equation}
Here \(\tau_k\) denotes the level-(k) local analyticity time scale starting from (t); that is, the length of the forward time interval on which the complexified Duhamel estimate for \(D^k u\) is valid. With this notation, the time-weighted bridge inequality can be written as
\begin{equation}\label{time-weight-Mjk}
\norm{D^j u(t)}{}^{\frac{1}{j+1}}
\leq
\mathcal{M}_{j,k}
\norm{D^k u(t)}{}^{\frac{1}{k+1}},
\qquad
\ell\leq j\leq k.
\end{equation}

The resulting analyticity estimate is stated in Corollary \ref{Analytic-estimate}.

\begin{remark}
   Since the above time-weighted bridge inequalities allow for more singular behavior, the time of local existence and -- in turn -- the lower bound on the radius
   of spatial analyticity will be weaker than in the case when the ratios were bounded by constants. However, it will turn out that it is still strong enough to dominate 
   the corresponding scale of sparseness for sufficiently large $k$.
\end{remark}

\begin{remark}
	The utility of Corollary \ref{Analytic-estimate} (as was the case for its counterpart in \cite{farhat}) is in applying it either along a sequence of times approaching the possiple singular time or at a single escape time close enough to the possible singular time. Since our application in the proof of the subsequent theorem would be of the latter type, we state it in terms of a single time increment.
\end{remark}

\begin{corollary}
\label{Analytic-estimate}
    Let \(\beta>1, \mu > 1, u_0 \in L^2, \)  and \(D^iu_0 \in L^{\infty}\) for \(0 \le i \le k\). Suppose that
\begin{align}\label{bridge}
    \norm{D^ju_0}{}^{\frac{1}{j+1}} \le (\tau_k)^{-\tfrac{1}{2\beta} \left(\tfrac{\mu + j}{j+ 1} - \tfrac{\mu + k}{k + 1}\right)} \dfrac{\left(j!\right)^{\frac{1}{j+1}}}{\left(k!\right)^{\tfrac{1}{k+1}}}  \norm{D^ku_0}{}^{\tfrac{1}{k+1}}, \qquad \ell \le j \le k,
\end{align}
where \(\ell\) and \(k\) are such that
\begin{align}\label{l!}
    \ell! \le \sqrt{\norm{u_0}{}} \le \left(k!\right)^{\tfrac{1}{k+1}},
\end{align}
and the \(j\)-dependent threshold $\tau_j$ is given by
\[
\tau_j:=
\left(\left[1 + \norm{u_0}{2}
\left(1 +\dfrac{1}{2} \log \norm{u_0}{}\right)\right]
\norm{D^ku_0}{} (j!)^{k - j}\right)^{-\tfrac{1}{k+1} \tfrac{2\beta(k+1)}{(2\beta - 1 )(k+1)- (\mu -1) 2k}}.
\]
Then the complexified solution satisfies 
\begin{align*}
    \sup_{t\in(0,\tau_j)} \;\sup_{y\in \mathcal{D}_t} \norm{D^{j}u(\cdot,y,t)}
  \;+\;
  \sup_{t\in(0,\tau_j)} \;\sup_{y\in \mathcal{D}_t} \norm{D^{j}v(\cdot,y,t)}
  \;\lesssim\;
  \norm{D^{j}u_0}{}  +  \\ \delta_0 \dfrac{\left(j!\right)^{\tfrac{1}{j+1}}}{\left(k!\right)^{\tfrac{1}{k+1}}} \norm{D^{k}u_0}{}^{\tfrac{j+1}{k+1}}(\tau_k)^{1-\tfrac{1}{2\beta} - \tfrac{\mu-1}{2\beta}\tfrac{2k}{k+1}},
\end{align*}
with the region of analyticity \(\Om_t\) defined by
\[\Om_t := \left\{ z = x+iy \in \mathbb{C}^3 : |y|\le ct^{\frac{1}{2\beta}}\right\}.\]
\begin{proof}
    It suffices to check conditions \eqref{cond1} and \eqref{cond2} as in Theorem~\ref{theorem1}.   
For the first condition \eqref{cond1}, we have  for all \(2\ell \le j \le k\)
\begin{align*}
    \sum_{\ell = i}^{ j-\ell}
  \binom{j}{i}\,
  \M_{i,k}^{\,i+1}\,
  \M_{j-i,k}^{\,j-i+1} &= \sum_{\ell = i}^{ j-\ell}
  \dfrac{j!}{(j-i)!i!} 
  \left[(\tau_k)^{-\tfrac{\mu-1}{2\beta}\left(\tfrac{1}{i+1} - \tfrac{1}{k+1} \right)} \dfrac{(i!)^{\tfrac{1}{i+1}}}{(k!)^{\tfrac{1}{k+1}}}\right]^{i+1}  \times \\ &
  \: \left[(\tau_k)^{-\tfrac{\mu-1}{2\beta}\left(\tfrac{1}{j-i+1} - \tfrac{1}{k+1} \right)} \dfrac{((j-i)!)^{\tfrac{1}{i+1}}}{(k!)^{\tfrac{1}{k+1}}}\right]^{j-i+1}\\
  &\le \dfrac{jj!}{(k!)^{\tfrac{j+2}{k+1}}} (\tau_k)^{-\tfrac{\mu-1}{2\beta} \tfrac{2k -j}{k+1}} \\
  & \lesssim\; \dfrac{(j!)^{\tfrac{1}{j+1}}}{(k!)^{\tfrac{1}{k+1}}} (j!)^{\tfrac{k-j}{k+1}} (\tau_k)^{ - \tfrac{\mu - 1}{2\beta} \tfrac{2k}{k+1}},
\end{align*}
The last inequality is valid due to the following:
Since \(j \le k\), then \(2k-j \in [k, 2k]\). Then for \((\tau_k) \in (0,1)\), we have 
\[(\tau_k)^{-\tfrac{\mu - 1}{2\beta} \tfrac{k}{k+1}} \le (\tau_k)^{-\tfrac{\mu - 1}{2\beta} \tfrac{2k-j}{k+1}} \le (\tau_k)^{-\tfrac{\mu - 1}{2\beta} \tfrac{2k}{k+1}}.\]
Also, \(l \le  j \le k\), so
\begin{align*}
    \dfrac{jj!}{(k!)^{\tfrac{j+2}{k+1}}} = \dfrac{jj!}{(k!)^{\tfrac{1}{k+1}} (k!)^{\tfrac{j+1}{k+1}}} \le \dfrac{jj!(j!)^{-\tfrac{j+1}{k+1}}}{(k!)^{\tfrac{1}{k+1}}} \le  \dfrac{(j!)^{\tfrac{1}{j+1}}}{(k!)^{\tfrac{1}{k+1}}} (j!)^{\tfrac{k-j}{k+1}}.
\end{align*}

Similarly, to verify that condition \eqref{cond2} holds, we compute
\begin{align*}
    \norm{u_0}{2} \!
  \sum_{0 \le i \le \ell}
  \binom{j}{i}
  \M_{\ell,k}^{\,\frac{(\ell+1)(i+\frac{3}{2})}{\ell+\frac{3}{2}}}
  \M_{j-i,k}^{\,j-i+1}
  \Big((k!)^{\frac{1}{k+1}}\norm{u_0}{}\Big)^{\frac{\left(\frac{3}{2}-1\right)(\ell-i)}{\ell+\frac{3}{2}}} \\ \le
  \norm{u_0}{2} \sum_{i = 0}^{l} \dfrac{j!}{(j-i)!i!} 
  \left[(\tau_k)^{-\frac{\mu-1}{2\beta}\left(\frac{1}{i+1} - \frac{1}{k+1} \right)} \dfrac{(i!)^{\frac{1}{i+1}}}{(k!)^{\frac{1}{k+1}}}\right]^{\frac{(l+1)(i+3/2)}{l+3/2}}  \times \\
   \left[(\tau_k)^{-\frac{\mu-1}{2\beta}\left(\frac{1}{j-i+1} - \frac{1}{k+1} \right)} \dfrac{((j-i)!)^{\frac{1}{i+1}}}{(k!)^{\frac{1}{k+1}}}\right]^{j-i+1} \left((k!)^{\frac{1}{k+1}} \norm{u_0}{}\right)^{\frac{(3/2-1)(l-i)}{l+3/2}}\\
  = \norm{u_0}{2} \sum_{i=0}^{l} \dfrac{j!(l!)^{\frac{i+3/2}{l+3/2}}}{i!} 
  \left( \dfrac{1}{(k!)^{\frac{1}{k+1}}}\right)^{(l+1)\frac{i+3/2}{l+3/2} + j-i+1} 
  \left((k!)^{\frac{1}{k+1}} \cdot \norm{u_0}{}\right)^{\frac{1}{2}\left(1 - \frac{i+3/2}{l+3/2}\right)} %^{\frac{1}{2}\left(1 - \frac{i+3/2}{l+3/2}\right)}
  (\tau_k)^{**},
\end{align*}
where the exponent \((**) \text{ is}  -\tfrac{\mu - 1}{2\beta} \left(\tfrac{(k-l)(i+3/2)}{l+3/2}  + k-j+i\right)\tfrac{1}{k+1}  \).

Applying \(l! \le \sqrt{\norm{u_0}{}} \le (k!)^{\tfrac{1}{k+1}}\), we have
\begin{align*}
    LHS &\lesssim \norm{u_0}{2} \dfrac{j! \sqrt{\norm{u_0}{}}}{(k!)^{\tfrac{j+2}{k+1}}} \sum_{i=0}^{l} \dfrac{1}{i!} \left( \dfrac{l!}{\sqrt{\norm{u_0}{}}}\right)^{\tfrac{i+3/2}{l+3/2}} (\tau_k)^{-\tfrac{\mu - 1}{2\beta} \left(\tfrac{(k-l)(i+3/2)}{l+3/2}  + k-j+i\right)\tfrac{1}{k+1} }\\
    &\lesssim \norm{u_0}{2} \dfrac{(j!)^{\tfrac{1}{j+1}}}{(k!)^{\tfrac{1}{k+1}}} (j!)^{\tfrac{k-j}{k+1}} (l+1)(\tau_k)^{-\tfrac{\mu - 1}{2\beta} \tfrac{2k - j}{k + 1}}\\
    &\lesssim  \norm{u_0}{2} \left(1 +\dfrac{1}{2} \log \norm{u_0}{}\right)\dfrac{(j!)^{\tfrac{1}{j+1}}}{(k!)^{\tfrac{1}{k+1}}} (j!)^{\tfrac{k-j}{k+1}} (\tau_k)^{-\tfrac{\mu - 1}{2\beta} \tfrac{2k}{k + 1}}.
\end{align*}
The last inequality holds for sufficiently large \(l\), in particular for \( l \geq  e^2\).
The Stirling's formula yields
\begin{align*}
    l! \sim  \sqrt{2\pi l } \left(\dfrac{l}{e}\right)^{l} \geq  e^l\sqrt{2\pi l} \implies  e^l \lesssim l!.
\end{align*}
Applying (\ref{l!}), we have
\begin{align*}
    e^l &\lesssim \sqrt{\norm{u_0}{}}.
   %l & \lesssim \dfrac{1}{2} \log\left( \norm{u_0}{} \right)
\end{align*}

Therefore, \( l  \lesssim \dfrac{1}{2} \log\left( \norm{u_0}{} \right)\).

Let \[\phi(j,k) = \dfrac{(j!)^{\frac{1}{j+1}}}{(k!)^{\frac{1}{k+1}}} (j!)^{\frac{k-j}{k+1}}\:\: \text{ and  }\:\: \psi(j,k) =\norm{u_0}{2}
\left(1 +\dfrac{1}{2} \log \norm{u_0}{}\right) \phi(j,k) .\]

 By Theorem~\ref{theorem1}, we have
    \begin{align*}
        \norm{D^ju(t)}{} \le \norm{D^ju_{0}}{} + \left(\tau_k\right)^{\tfrac{2\beta-1}{2\beta} - \tfrac{\mu - 1}{2\beta} \tfrac{2k}{k+1}} \left( \phi(j,k) + \psi(j,k)\right) \norm{D^ju_0}{}^{\tfrac{j+2}{k+1}}.
    \end{align*}
    
    We choose  $\tau_j$ such that
    \begin{align*}
        (\tau_j)^{\tfrac{(2\beta-1)(k+1) - (\mu -1)2k}{2\beta (k+1)}} \le 
        \left( 1 + \norm{u_0}{2}\left(1 +\dfrac{1}{2} \log \norm{u_0}{}\right) \right)^{-1} 
        \norm{D^ku_0}{}^{- \tfrac{1}{k+1}}
        \left( j!\right)^{-\tfrac{k-j}{k+1}}
    \end{align*}
 completing the proof.
\end{proof}
\end{corollary}
We now compare the radius of analyticity with the scale of sparseness. For the radius of  analyticity, we take the  
\(2 \beta\)-th root of the existence time
\begin{align}\label{rhok}
    \rho_k \geq \dfrac{1}{c(\norm{u_0}{2}, \log(\norm{u_0}{}), \beta, \mu)} \dfrac{1}{\norm{D^k u_0}{}^{\tfrac{1}{2\beta -1 - (\mu -1)\cdot \tfrac{2k}{k+1}}\cdot \tfrac{1}{k+1}}}.
\end{align}

On the other hand, Theorem \ref{scalebound} yields the following lower bound for the sparseness scale:
\[
r_k \ge \frac{1}{\|D^k u_0\|^{\tfrac{1}{k+3/2}}}.
\]
The fundamental requirement is that \(\rho_k \geq r_k \), that is, the analyticity radius dominates the sparseness scale.  It suffices to take \(2\beta - 1 -(\mu - 1) \cdot \dfrac{2k}{k+1} > 1\), which holds for \(\mu < 2\beta -1\). In particular, the threshold advances from \(\mu<2\beta-1\) to \(\mu<\beta\) (see proof of Theorem \ref{main-theorem}).

With the analyticity estimate and the analyticity--sparseness comparison established, we can now complete the proof of the main theorem.

\subsection{Proof of Theorem \ref{main-theorem}}

\begin{proof}
    Let \(k \in \Z^{+}\) and \(0 \le j \le k\). By Taylor's theorem, we have \(c_{\alpha, k}^{i}(t) = \dfrac{1}{\alpha!} \dfrac{\partial^k }{\partial^\alpha x} u^{i}(x^{*}, t) \). Applying condition (ii) and taking the maximizing component \((i, \alpha)\),  
    we have
    \begin{align*}
        \dfrac{\|D^{j}u(t) \|^{\tfrac{1}{j+1}}}{\|D^{k}u(t) \|^{\tfrac{1}{k+1}}} = 
         \dfrac{|D^{j}u(t) |^{\tfrac{1}{j+1}}}{|D^{k}u(t) |^{\tfrac{1}{k+1}}} \le 
         \dfrac{(j!)^{\tfrac{1}{j+1}}}{(k!)^{\tfrac{1}{k+1}}} \dfrac{(c_{\alpha, j}^{i})^{\tfrac{1}{j+1}}}{(c_{\alpha, k}^{i})^{\tfrac{1}{k+1}}}.
    \end{align*}
    By (i), for a constant $C > 1$ uniform in $i$ and $|\alpha| =k$ we have 
    \begin{align*}
         \dfrac{\|D^{j}u(t) \|^{\tfrac{1}{j+1}}}{\|D^{k}u(t) \|^{\tfrac{1}{k+1}}} 
         &\le  \dfrac{(j!)^{\tfrac{1}{j+1}}}{(k!)^{\tfrac{1}{k+1}}} \dfrac{\left[\delta_{\alpha, j}^{i}(t) \rho(t)^{-j}\right]^{\tfrac{1}{j+1}}}{\left[\delta_{\alpha, k}^{i}(t) \rho(t)^{-k}\right]^{\tfrac{1}{k+1}}}\\
         & \le C^{\tfrac{j}{j+1} +\tfrac{k}{k+1}} 
          \rho(t)^{\tfrac{k}{k+1} - \tfrac{j}{j+1}}(T^{*} - t)^{ - \tfrac{\gamma }{\log j} 
          \tfrac{j}{j+1} 
          - \tfrac{\gamma }{\log k} 
          \tfrac{k}{k+1} }
          \dfrac{(j!)^{\tfrac{1}{j+1}}}{(k!)^{\tfrac{1}{k+1}}}\\
          & \le  C^{\tfrac{j}{j+1} +\tfrac{k}{k+1}} 
          \rho(t)^{\tfrac{k}{k+1} - \tfrac{j}{j+1}}
          (T^{*} -t)^{-\tfrac{2\gamma}{\log k}}
          \dfrac{(j!)^{\tfrac{1}{j+1}}}{(k!)^{\tfrac{1}{k+1}}}.
    \end{align*}
Since  \(\rho \to 0\) as \(t \to T^{*}\), we choose \(t_0 \in (0, T^{*})\) so that \(\rho(t) \le 1\) for all \(t \in (t_0, T^{*}) \). Then  \(\rho(t)^{\tfrac{k}{k+1} - \tfrac{j}{j+1}}
\le 1 \) for all \(t \in (t_0, T^{*})\). 
Also, \(C^{\tfrac{j}{j+1} + \tfrac{k}{k+1}} \le C^2\). Choosing \(\gamma \sim \tfrac{\mu - 1}{4\beta}\) and restricting to \( j \in [\log k,  k] \), with \(j\) and \(k\) satisfying \eqref{l!} one has
\begin{align*}
    \dfrac{\|D^{j}u(t) \|^{\tfrac{1}{j+1}}}{\|D^{k}u(t) \|^{\tfrac{1}{k+1}}}  
    \lesssim 
    (T^{*} -t)^{-\tfrac{\mu-1}{2\beta} \cdot\tfrac{k - j}{(j+1)(k + 1)}}
    \dfrac{(j!)^{\tfrac{1}{j+1}}}{(k!)^{\tfrac{1}{k+1}}},
\end{align*}
which is compatible with \eqref{bridge}. (since we are in a blow-up scenario, \(T^{*}-t > \tau_k \)).

Let \(k\) be a fixed constant and let \(t \in (T^{*}- \varepsilon, T^{*})\) be an escape time for \(D^ku\).  That is to say, 
\begin{align}\label{escapetime}
    \norm{D^{k}u(\tau)}{} > \norm{D^ku(t)}{}\:\:   \text{  for all } \:\: \tau \in (t, T^{*}).
\end{align}

To analyze the case for \(j = k\), we evolve the system from time \(t\) to \(s= t + \tau_{k}.\)  First, observe that the following holds:
\[1 + \norm{u_0}{2}(1 + \tfrac{1}{2}\log(\norm{u(t)}{})) \leq (1 + \norm{u_0}{2})(1 + \log\norm{u(t)}{} ).\]
Using this inequality, we have the following for the time increment $\tau_k$
\begin{align*}
    \tau_k &= \left(\left[1 + \norm{u_0}{2}
    \left(1 +\dfrac{1}{2} \log \norm{u(t)}{}\right)\right]
    \norm{D^ku(t)}{}\right)^{-\tfrac{1}{k+1} \tfrac{2\beta(k+1)}{(2\beta - 1 )(k+1)- (\mu -1) 2k}}\\
    &\geq \dfrac{1}{\tilde{c}(\norm{u_0}{2}, \beta, \mu, \delta_0)} 
    \left(\log(\norm{u(t)}{}) \norm{D^ku(t)}{}\right)^{-\tfrac{2\beta(k+1)}{(2\beta - 1 )(k+1)- 2(\mu -1)k} \tfrac{1}{k+1}}.
\end{align*}

Similarly, the upper bound on the complexified solution obtained in Corollary~\ref{Analytic-estimate} yields

\begin{align}\label{supbound}
    \sup_{y\in \mathcal{D}_s} \norm{D^{k}u(\cdot,y,s)}
  \;+\;
   \sup_{y\in \mathcal{D}_s} \norm{D^{k}v(\cdot,y,s)}
  \;\lesssim\; (1  + \delta_0 (\tau_k)^{1-\tfrac{1}{2\beta} - \tfrac{\mu-1}{2\beta}\tfrac{2k}{k+1}})  \norm{D^{k}u(t)}{}.
\end{align}

Using the harmonic measure maximum principle we can prevent the blow-up. This is done by making sure that the radius of spatial analyticity at \(s,  \rho_k\) dominates the a priori scale of sparseness at \(s, r_k\) \cite{Bradshaw}.  By \eqref{rhok} and Theorem \ref{scalebound}, we obtain the radius of analyticity \(\rho_k\) and the sparseness scale \(r_k\), respectively, as
\begin{align*}
	\rho_k &= \dfrac{1}{c_2(\norm{u_0}{2}, \beta, \mu, \delta_{0})}
	\left(\log(\norm{u(s)}{})\norm{D^{(k)}u(s)}{}\right)^{-a_k},
	\\
	r_k &= c_1 (\norm{u_0}{2}) \norm{D^{(k)}u(s)}{}^{-\tfrac{1}{k+ \tfrac{3}{2}}}.
\end{align*}

Where \(a_k: = \dfrac{1}{2\beta -1 - (\mu -1) \tfrac{2k}{k+1} } \dfrac{1}{k+1}\).

Observe that the difference in exponents  
\begin{align}\label{gap}
    \dfrac{1}{k+\tfrac{3}{2}} - a_k 
    %{\dfrac{1}{2\beta -1 - (\mu -1) \tfrac{2k}{k+1} } \dfrac{1}{k+1}} 
    =
\dfrac{(2\beta - 2\mu)k +(2\beta - 5/2)}{(k+3/2) [(k+1)(2\beta -1)- (\mu -1)2k]} =: \Delta_k
\end{align}

is positive for all \(k > \dfrac{5/2 - 2\beta}{2\beta -2\mu} \) for all \(\beta \in (1, 5/4)\). This holds because of the following argument: For \(\beta \in (1, 5/4)\), \(2\beta - 5/2 < 0\). Therefore, the numerator is positive for all  \(\mu < \beta\). Therefore, 
\( k > \dfrac{5/2 - 2\beta}{2\beta - 2\mu}\) does the job. In the denominator, for large \(k\)
\[(k+1) (2\beta - 1)- (\mu -1 )2k \sim k(2\beta- 1) - 2(\mu -1)k,\]
which is positive if \(\mu < \beta + \dfrac{1}{2}\). This sufficient condition on \(\mu\) is automatically satisfied since \(\mu < \beta\). Consequently,
\begin{align}\label{dominates}
  \rho_k(s)
  \;\ge\; (c_1c_2)^{-1}(\log\|u(s)\|)^{-a_k}\,\|D^{k}u(s)\|^{\Delta_k}\, r_k(s)
  \;\ge\; r_k(s).
\end{align}
The last inequality holds if
\[
\|D^{k}u(s)\|^{\Delta_k}\ \ge\ c_1c_2\,(\log\|u(s)\|)^{a_k}.
\]
Since $a_k>0$ and $\log\|u(s)\|\ge1$ near $T^*$, we can choose $k$  large enough so that this is satisfied. Hence, the analytic tube contains every complex disk with radius $\rho\le r_k(s)$.

Let \(\xi \in \mathbb{S}^2\) be the direction from the \(1D \: \delta\)-sparseness at the time \(s\) and consider the complex line \(z \mapsto x_0 + z\xi.\) Since \(\rho \leq r_k(s) \leq \rho_k(s),\) the map is contained in the analyticity tube at time \(s.\) 

Define the scalar analytic function
\[f(z): = (D^{(k)}u)_{j}^{\pm}(x_0 + z \xi, s), \qquad z \in \mathbb{D}_{\rho}.\]

From the tube bound (obtained in the analyticity step (\ref{supbound})), there exists \(M\) such that

\begin{equation}\label{eq:supM-3d}
\sup_{|z|<\rho}|f(z)|\ \le\ M,\qquad M:=A_k\,\|D^k u(t)\|,
\end{equation}
where $A_k:=1+\delta_0 (\tau_k)^{1-\tfrac{1}{2\beta}-\tfrac{\mu-1}{2\beta}\tfrac{2k}{k+1}}$.
On the real diameter 
$I = [-\rho,\rho]\subset\partial\mathbb{D}_\rho$, the $1$D $\delta$-sparseness
and the definition of  $V^{j,\pm}_\lambda (s)$ imply that on a “good” subset of length at least
$(1-\delta)2\rho$ (i.e., \(|I \setminus V_{\lambda}^{j,\pm}(s)| \geq (1-\delta )2\rho\)),
\[
|f|\ \le\ \lambda\,\|D^k u(s)\|\ \le\ \lambda\,M,
\]
while everywhere $|f|\le M$ by \eqref{eq:supM-3d}. 

Next, we apply the harmonic–measure maximum
principle to $\log|f|$ in the slit disk
\(\mathbb{D}_\rho\setminus\big([-\,\rho,\rho]\cap V^{j,\pm}_\lambda(s)\big)\). Let \(h = h(\delta)\)  be the  harmonic measure of the good arc at the center. Then 
\begin{align*}
    \log|f(0)| & \leq h \log (\lambda M) + (1-h) \log(M) = \log(\lambda^hM),\\
   \implies  |f(0)| &\leq \lambda^h M.
\end{align*}

Now, pick the maximizing component/sign and \(x_0\) as in the statement of the theorem. Putting \(|f(0)| = \norm{D^{k}u(s)}{}\) and  \(M = A_k \norm{D^{k}u(t)}{} \), we get
\[  \norm{D^{k}u(s)}{} \leq \eta  \norm{D^{k}u(t)}{}, \qquad \text{where } \eta: = \lambda^h A_{k}.\]

Observe that \(A_k \to 1\) as \(t\) approaches \(T^*\). Take \(t\) close enough to \(T^*\) such that \(A_k \leq 1  + \varepsilon\) and choose \((\lambda,  \delta)\) such that
\[\lambda^h (1 + \varepsilon) < 1 \implies \eta < 1.\]

Therefore, we obtain the strict contraction
\begin{equation}\label{strict-contraction}
\|D^k u(s)\|_{\infty}
<
\|D^k u(t)\|_{\infty}.
\tag{21}
\end{equation}
On the other hand, by the assumed monotonicity of the coefficient functions
\(c_{\alpha,k}^{i}\), together with the focused-buildup condition, every
\(t\in(T^*-\varepsilon,T^*)\) satisfying 
\[\|u(t)\|\leq k^2
\]
(arising from \eqref{cond2}) has the escape-time property at level \(k\). Consequently,
\[
\|D^k u(s)\|_{\infty}
>
\|D^k u(t)\|_{\infty}
\qquad
\text{for every }s\in(t,T^*),
\]
and in particular, for the time \(s\in(t,t+\tau_k)\) considered above. This contradicts \eqref{strict-contraction}.

Since \(t\in(T^*-\varepsilon,T^*)\) was arbitrary and the argument applies at every sufficiently large differential level \(k\) for which \(\Delta_k>0\), the analyticity radius remains positive up to \(T^*\). Hence, the solution extends analytically beyond \(T^*\), and therefore \(T^*\) cannot be a blow-up time.
\end{proof}

\section*{Acknowledgments}

The author wishes to thank Professor Z. Grujić for proposing the problem, for his continued support, and for many inspiring conversations.

%\vspace{1cm}
%The University of Alabama at Birmingham\\
%\href{mailto:mpphiri@uab.edu}{mpphiri@uab.edu}
\end{document}